\documentclass[
aps,%
12pt,%
final,%
notitlepage,%
oneside,%
onecolumn,%
nobibnotes,%
nofootinbib,% 
superscriptaddress,%
noshowpacs,%
centertags]%
{revtex4}

\usepackage{amsmath,amssymb}
\usepackage{enumerate}
\usepackage{graphicx}
\usepackage{xcolor}
\usepackage{mdframed}

\newmdenv[
linewidth=2pt,
linecolor=gray,
  topline=false,
  bottomline=false,
  rightline=false,
  skipabove=\topsep,
  skipbelow=\topsep
]{siderules}

\newcommand{\R}{\mathbb{R}}

\renewcommand{\P}{\operatorname{\mathsf{P}}}
\newcommand{\Var}{\operatorname{\mathsf{Var}}}
\newcommand{\Cov}{\operatorname{\mathsf{Cov}}}
\newcommand{\1}{%\operatorname
{\mathsf{1}}}

\newcommand{\even}{\mathsf{even}}
\newcommand{\odd}{\mathsf{odd}}

\newcommand{\ICCS}{\text{ICCS}_n}

\newcommand{\Si}{\Sigma}

\begin{document}
%\selectlanguage{english}
\begin{center}
\large \bfseries %\MakeTextUppercase
{%
%\color{blue}
What Intraclass Covariance Structures Can Symmetric Bernoulli Random Variables Have?
}
\end{center}
\begin{center}
%\bfseries 
Iosif Pinelis
\end{center}

%\vspace*{.1cm}

\begin{center}
	{\it Michigan Technological University, Houghton, Michigan, USA}
	
	 Received June 16, 2022 \\ 
Revised   September 24, 2022  \\ 
Accepted October 3, 2022

\end{center}

\begin{center}
\begin{minipage}{5.5in}
\small
{\bf%\color{blue} 
Abstract}---The covariance matrix of random variables $X_1,\dots,X_n$ is said to have an intraclass covariance structure if the variances of all the $X_i$'s are the same and all the pairwise covariances of the $X_i$'s are the same. We provide a possibly surprising characterization of such covariance matrices in the case when the $X_i$'s are symmetric Bernoulli random variables. 

\medskip

{%\color{blue}
} 
\end{minipage}
\end{center}

\vspace*{.6cm}

For natural $n\ge2$, let $\Si=[\Si_{i,j}]_{i,j\in[n]}$ be the covariance matrix of random variables (r.v.'s) $X_1,\dots,X_n$ with finite second moments, so that $\Si_{i,j}=\Cov(X_i,X_j)$ for all $i$ and $j$ in the set $[n]:=\{1,\dots,n\}$. We are assuming that the matrix $\Si$ is nonzero.  

The covariance matrix $\Si$ is said to have an intraclass covariance structure if (i) $\Si_{i,i}=\Var X_i=
\Cov(X_i,X_i)$ is the same for all $i\in[n]$ and (ii) $\Si_{i,j}=\Cov(X_i,X_j)$ is the same for all distinct $i$ and $j$ in $[n]$. Let $\ICCS$ denote the set of all $n\times n$ covariance matrices that have an intraclass covariance structure. 

In particular, if the r.v.'s $X_1,\dots,X_n$ are exchangeable -- that is, if the joint distribution of the $X_i$'s is invariant with respect to all permutations of the indices $1,\dots,n$ (see e.g. \cite{kallenberg-exch} for much more on exchangeability of r.v.'s), then the covariance matrix $\Si$ will be in the set $\ICCS$. So, one may say that the covariance matrix $\Si$ has 
an intraclass covariance structure if the r.v.'s $X_1,\dots,X_n$ pertain to items that belong to one class and thus are exchangeable in a certain weak sense; this explains the use of the term ``intraclass''. 
The notion of an intraclass covariance structure was introduced by Fisher \cite{fisher-book} and has been studied in many subsequent papers, including e.g. \cite{walsh47,press,srivastava}. 

Obviously, the covariance matrix $\Si$ is in the set $\ICCS$ if and only if 
\begin{equation}
	\Si=(a-b)I_n+b\,\1_n\1_n^\top \label{eq:Si=}
\end{equation}
for some real numbers $a$ and $b$, where $I_n$ is the $n\times n$ identity matrix and $\1_n:=[1,\dots,1]^\top$, the $n\times1$ matrix of $1$'s. 

Recall that a real $n\times n$ matrix is a covariance matrix if and only if it is positive semidefinite; cf.\ e.g.\ \cite[Sect. III.6, Theorem 4]{feller_vol2}. Note that (i) $\1_n$ is an eigenvector of the matrix $\1_n\1_n^\top$ belonging to the eigenvalue $n$ and (ii) any nonzero vector orthogonal to $\1_n$ is an eigenvector of the matrix $\1_n\1_n^\top$ belonging to the eigenvalue $0$. 
So, the only eigenvalues of the matrix $\Si$ of the form \eqref{eq:Si=} are $a-b+bn$ and $a-b$. 

It follows that the matrix $\Si$ of the form \eqref{eq:Si=} is in $\ICCS$ if and only if $-\frac a{n-1}\le b\le a$, that is, if and only if the pairwise correlation, $\rho=b/a$, between r.v.'s whose covariance matrix has an intraclass covariance structure is no less that $-1/(n-1)$: 
\begin{equation}
	\rho\ge\rho_{n,\min}:=-\frac1{n-1}.  \label{eq:rho ge}
\end{equation}
This is in contrast with the general lower bound $-1$ on the correlation between arbitrary r.v.'s. Let us refer to the values of $\rho$ satisfying condition \eqref{eq:rho ge} as \emph{good}. 

In the rest of this note, we shall consider the special case when the r.v.'s $X_1,\dots,X_n$ are symmetric Bernoulli, so that 
\begin{equation}
	\P(X_i=1)=\tfrac12=\P(X_i=0) \label{eq:bern}
\end{equation}
for all $i\in[n]$. This important case has been extensively studied in computer science in general and in machine learning in particular (see e.g. \cite{Nnatarajan13,senel18,%meng20,
baldi-vershynin21,zhang-etal21}), as well as in other applications of probability theory -- though mainly when the $X_i$'s are independent. 

%
%Condition \eqref{eq:bern} implies $\E X_i=\tfrac12$ and $\Var X_i=\tfrac14$ for all $i\in[n]$. Moreover, for all distinct $i$ and $j$ in $[n]$, 
%\begin{equation*}
%	\Cov(X_i,X_j)=\E X_iX_j-(\tfrac12)^2=\P(X_i=X_j=1)-\tfrac14  
%\end{equation*}
%and 
%\begin{equation*}
%\begin{aligned}
%\P(X_i=X_j=1)&=\P(X_i=1)-\P(X_i=X_j=0) \\ 
%&=\P(X_i=0)-\P(X_i=X_j=0)=\P(X_i=X_j=0), 
%\end{aligned}
%\end{equation*}
%so that $\P(X_i=X_j=1)=\P(X_i=X_j=0)=\tfrac12\,P(X_i=X_j)$ and hence 
%\begin{equation*}
%	\Cov(X_i,X_j)=\tfrac12\,P(X_i=X_j)-\tfrac14.   
%\end{equation*}
%Hence, given condition \eqref{eq:bern}, the probabilities $P(X_i=X_j)$ are the same for all distinct $i$ and $j$ in $[n]$, so that there is some $p\in[0,1]$ such that  
%\begin{equation}
%	\text{$P(X_i=X_j)=p$ for all distinct $i$ and $j$ in $[n]$.} \label{eq:p}
%\end{equation} 
%
%Therefore, given condition \eqref{eq:bern}, condition \eqref{eq:a,b} can be rewritten as 
%\begin{equation}
%	p\ge\frac{n-2}{2(n-1)}. \label{eq:p ge}
%\end{equation}
%That is, the covariance matrix of symmetric Bernoulli r.v.'s $X_1,\dots,X_n$ is in the set $\ICCS$ if and only if conditions \eqref{eq:p} and \eqref{eq:p ge} hold. 

%the probabilities $P(X_i=X_j)$ 
%are the same for all distinct $i$ and $j$ in $[n]$ and \eqref{eq:p ge} holds for $p$ defined by \eqref{eq:p}. 

The question now is the following: 

\bigskip

%\framebox{
%\begin{minipage}{5.5in}
\begin{siderules}
For what values of pairwise correlation $\rho$ do there exist symmetric Bernoulli r.v.'s $X_1,\dots,X_n$ whose covariance matrix $\Si$ is in $\ICCS$?
\end{siderules}
%\end{minipage}
%} 

\bigskip

Let us refer to such values of $\rho$ as \emph{symmetric-binary-good}. Clearly, any symmetric-binary-good value of $\rho$ must be good. One then may wonder  
%the question is then 
whether every good value of $\rho$ is symmetric-binary-good. 

The answer to this question may seem surprising: 
\begin{itemize}
	\item if $n$ is even, then yes, every good value of $\rho$ is symmetric-binary-good; 
	\item if $n$ is odd, then ``nearly every'' good value of $\rho$ is symmetric-binary-good.  
\end{itemize}

For symmetric Bernoulli r.v.'s $X_1,\dots,X_n$ whose covariance matrix $\Si$ is in $\ICCS$, it is a bit more convenient to deal with the probability 
\begin{equation*}
	p:=\P(X_1=X_2)
\end{equation*}
than with the correlation $\rho$. It is easy to see that the values of $\rho$ and $p$ are in the simple bijective correspondence 
\begin{equation}
	(-1,1)\ni 2p-1=\rho\longleftrightarrow p=\frac{1+\rho}2\in(0,1), \label{eq:p,rho}
\end{equation}
so that $\P(X_i=X_j)=p$ for all distinct $i$ and $j$ in $[n]$. 

Let us refer to the values of $p$ corresponding to the good values of $\rho$ as \emph{good} values of $p$, and let us similarly define the \emph{symmetric-binary-good} values of $p$. So, in view of \eqref{eq:rho ge} and \eqref{eq:p,rho}, a value $p\in(0,1)$ is good if and only if 
\begin{equation}
	p\ge p_n:=\frac{n-2}{2(n-1)}. \label{eq:p ge}
\end{equation}
Thus, we have to determine the symmetric-binary-good values of $p$. 

Suppose for a moment that $p\in(0,1)$ is symmetric-binary-good. Then there exist symmetric Bernoulli r.v.'s $X_1,\dots,X_n$ such that $\P(X_i=X_j)=p$ for all distinct $i$ and $j$ in $[n]$. Letting $g$ stand for the joint probability mass function of the r.v.'s $X_1,\dots,X_n$, we note that $g$ is a nonnegative function such that 
\begin{enumerate}[(i)]
	\item $\sum_{x\in\{0,1\}^n}g(x)=1$, 
\item $\sum_{x\in\{0,1\}^n}1(x_i=0)g(x)=\frac12$ for all $i\in[n]$, 
\item $\sum_{x\in\{0,1\}^n}1(x_i=x_j)g(x)=p$ for all distinct $i$ and $j$ in $[n]$;
\end{enumerate} 
of course, here $x_i$ denotes the $i$th coordinate of the vector $x=(x_1,\dots,x_n)\in\{0,1\}^n$. 
By symmetry, conditions (i)--(iii) will hold with $\tilde g(x):=\frac1{n!}\sum_{\pi\in \Pi_n}g(\pi(x))$ in place of $g(x)$, where $\Pi_n$ is the set of all permutations of the set $[n]$. Note that $\tilde g(x)=f(\sum_1^n x_i)$ for some nonnegative function $f\colon\{0,\dots,n\}\to\R$ and all $x\in\{0,1\}^n$. So, conditions (i)--(iii) can be rewritten as 
\begin{enumerate}[(I)]
	\item $\sum_{k=0}^n \binom nk f(k)=1$, 
\item  $\sum_{k=0}^n \binom{n-1}k f(k)=\frac12$ for all $i$, 
\item  $\sum_{k=0}^n a_{n,k} f(k)=p$,
\end{enumerate}
where 
\begin{equation*}
	a_{n,k}=\binom{n-2}k+\binom{n-2}{k-2}; 
\end{equation*}
of course, $\binom{n-1}n=0$, $\binom{n-2}k=0$ if $k\ge n-1$ and $\binom{n-2}{k-2}=0$ if $k\le1$. 

Thus, for any given $n\ge2$ and $p\in(0,1)$, we want to see whether there is a nonnegative function $f\colon\{0,\dots,n\}\to\R$ such that conditions (I)--(III) hold. 

Towards this goal, consider the problem of finding the extrema of $\sum_{k=0}^n a_{n,k} f(k)$ over all $f\in F_n$, where $F_n$ is the set of all nonnegative function $f\colon\{0,\dots,n\}\to\R$ satisfying condition (I). In view of the symmetries $\binom nk=\binom n{n-k}$ and $a_{n,k}=a_{n,n-k}$, without loss of generality the functions $f$ are symmetric in the same sense: $f(k)=f(n-k)$ for all $k\in\{0,\dots,n\}$---otherwise, replacing $f(k)$ by $\frac12\,(f(k)+f(n-k))$, %\break 
we will have the sums in (I) and (III) unchanged. 
Next, consider the ratios 
\begin{equation*}
r_k:=	r_{n,k}:=\frac{a_{n,k}}{\binom nk}
	=\frac{(n-k)(n-k-1)+k (k-1)}{n (n-1)}. 
\end{equation*}
Note that $r_{k+1}\le r_k$ if $0\le k\le\frac{n-1}2$ and $r_{k+1}\ge r_k$ if $\frac{n-1}2\le k\le n-1$. Also, $r_k=r_{n-k}$. So, 
%$r_0$ and $r_n$ are the largest ones among the $r_k$'s, and 
the smallest among the $r_k$'s is/are the one/ones with index/indices $k$ closest to $\frac n2$. 

%!!!! stopped here at 2:37 on 4/18/22 

More specifically, if $n=2m-1$ is odd, then $r_k\ge r_m=r_{m-1}$ for all $k\in\{1,\dots,n-1\}$. Letting then 
\begin{gather*}
	f^\odd_{\min}(m-1):=\frac{1/2}{\binom n{m-1}}=\frac{1/2}{\binom nm},\quad 
	f^\odd_{\min}(m):=\frac{1/2}{\binom nm}=\frac{1/2}{\binom n{m-1}}, \\ 
%\end{equation*}
%and 
%\begin{equation*}
	f^\odd_{\min}(k):=0\quad\text{for all}\quad k\in\{0,\dots,n\}\setminus\{m-1,m\}, 
\end{gather*}
we see that $f^\odd_{\min}$ is a symmetric function in $F_n$ and 
\begin{equation*}
	(r_k-r_m)(f^\odd_{\min}(k)-f(k))\le0
\end{equation*}
for all $k\in\{0,\dots,n\}$ and all symmetric functions $f\in F_n$,  
% -- because for such functions $f$ we have $f(0)=f(n)\le\frac12=f^\odd_{\min}(0)=f^\odd_{\min}(n)$. 
which implies  
\begin{equation*}
\begin{aligned}
	\sum_{k=0}^n a_{n,k} f^\odd_{\min}(k)-\sum_{k=0}^n a_{n,k} f(k)
	&=\sum_{k=0}^n a_{n,k} (f^\odd_{\min}(k)-f(k)) \\ 
	&=\sum_{k=0}^n \binom nk r_k (f^\odd_{\min}(k)-f(k)) \\ 
	&=\sum_{k=0}^n \binom nk (r_k-r_m)(f^\odd_{\min}(k)-f(k))\le0.  
\end{aligned}
\end{equation*}
It follows that $f^\odd_{\min}$ is a minimizer of $\sum_{k=0}^n a_{n,k} f(k)$ over all $f\in F_n$, that is, over all nonnegative $f$ satisfying condition (I). Moreover, condition (II) is satisfied with $f^\odd_{\min}$ in place of $f$. 

We conclude that, in the case when $n=2m-1$ is odd, 
$f^\odd_{\min}$ is a minimizer of $\sum_{k=0}^n a_{n,k} f(k)$ over all nonnegative $f$ satisfying \emph{both} conditions (I) and (II). The corresponding minimum value of $\sum_{k=0}^n a_{n,k} f(k)$ is 
\begin{equation*}
	p^\odd_{n,\min}:=\sum_{k=0}^n a_{n,k} f^\odd_{\min}(k)=\frac{m-1}{2m-1}=\frac{n-1}{2n}. 
\end{equation*}

Similarly, in the case when $n=2m$ is even, a minimizer of $\sum_{k=0}^n a_{n,k} f(k)$ over all nonnegative $f$ satisfying both conditions (I) and (II) is given by 
\begin{equation*}
	f^\even_{\min}(m):=\frac1{\binom nm}
%\end{equation*}
\quad\text{and}\quad  
%\begin{equation*}
	f^\even_{\min}(k):=0\ \ \,\text{for all}\ \,  k\in\{0,\dots,n\}\setminus\{m\}, 
\end{equation*}
and the corresponding minimum value of $\sum_{k=0}^n a_{n,k} f(k)$ is 
\begin{equation*}
	p^\even_{n,\min}:=\sum_{k=0}^n a_{n,k} f^\even_{\min}(k)=\frac{m-1}{2m-1}=\frac{n-2}{2(n-1)}. 
\end{equation*}

The above minimization can of course be recognized as something similar to, or even a special case of, the Neyman--Pearson lemma \cite[part III]{neyman-pearson}. 

The just considered cases of odd and even $n$ can be summarized as follows: For 
\begin{equation*}
	m_n:=\lceil n/2\rceil,
\end{equation*}
let $f_{\min}$ be the symmetric function in $F_n$ such that 
$%\begin{equation*}
	\sum_{k\in\{m_n,n-m_n\}}f_{\min}(k)=1,
$ %\end{equation*}
so that $f(k)=0$ for $k\in\{0,\dots,n\}\setminus\{m_n,n-m_n\}$. 
Then $f_{\min}$ is a minimizer of $\sum_{k=0}^n a_{n,k} f(k)$ over all nonnegative $f$ satisfying conditions (I) and (II). The corresponding minimum value of $\sum_{k=0}^n a_{n,k} f(k)$ is 
\begin{equation*}
	p_{n,\min}:=\sum_{k=0}^n a_{n,k} f_{\min}(k)=\dfrac{m_n-1}{2m_n-1}. 
\end{equation*}

The extremal joint distribution of the binary r.v.'s $X_1,\dots,X_n$ corresponding to the minimizer $f_{\min}$ can be described as follows: the random set $I:=\{i\in[n]\colon X_i=1\}$ is uniformly distributed on the set $S_n:=\binom{[n]}{m_n}\cup\binom{[n]}{n-m_n}$, where $\binom{[n]}k$ denotes the set of all subsets of cardinality $k$ of the set $[n]$; of course, $S_n:=\binom{[n]}{n/2}$ if $n$ is even. 

%of all subsets of the set $[n]$ whose cardinalities are in the set $\{m_n,n-m_n\}$. 

%the cardinality of the random set $\{i\in[n]\colon X_=1\}$ is in the set $\{m_n,n-m_n\}$ with probability $1$ and this random set is uniformly distributed on the set of all subsets of the set $[n]$ whose cardinalities are 
%\begin{equation*}
%	\{i\in[n]\colon X_=1\}
%\end{equation*}

Next, letting 
\begin{equation*}
	f_{\max}(0):=\tfrac12,\quad f_{\max}(n):=\tfrac12,\quad	f_{\max}(k):=0\ \, \text{for all}\ \,  k\in\{1,\dots,n-1\}, 
\end{equation*}
we see that the nonnegative function $f_{\max}$ satisfies conditions (I) and (II), and also 
$\sum_{k=0}^n a_{n,k} f_{\max}(k)=1$. On the other hand, for any nonnegative function $f$ satisfying conditions (I) and (II), the sum $\sum_{k=0}^n a_{n,k} f(k)$ is a probability and hence does not exceed $1$. 
We conclude that $f_{\max}$ is a maximizer of $\sum_{k=0}^n a_{n,k} f(k)$ over all nonnegative $f$ satisfying conditions (I) and (II). The corresponding maximum value of $\sum_{k=0}^n a_{n,k} f(k)$ is 
\begin{equation*}
	p_{n,\max}:=\sum_{k=0}^n a_{n,k} f_{\max}(k)=1. 
\end{equation*}

The extremal joint distribution of the binary r.v.'s $X_1,\dots,X_n$ corresponding to the maximizer $f_{\max}$ can be described as follows: the random set $I=\{i\in[n]\colon X_i=1\}$ is uniformly distributed on the set $\{\emptyset,[n]\}$; that is, $\P(I=\emptyset)=\frac12=\P(I=[n])$. 

Now note that the set of all values of $\sum_{k=0}^n a_{n,k} f(k)$, where $f\colon\{0,\dots,n\}\to\R$ is a nonnegative function such that conditions (I) and (II) hold, is convex and therefore coincides with the interval $[p_{n,\min},p_{n,\max}]=[p_{n,\min},1]$. 

Thus, a value $p\in(0,1)$ is symmetric-binary-good if and  only if 
\begin{equation*}
	p\ge p_{n,\min}=\dfrac{m_n-1}{2m_n-1}
	=\left\{
	\begin{alignedat}{2}
	&\dfrac{n-2}{2(n-1)}=p_n&&\text{\; if $n$ is even}, \\ 
	&\dfrac{n-1}{2n}=p_{n+1}>p_n&&\text{\; if $n$ is odd}, \\ 
	\end{alignedat}
	\right.
\end{equation*}
where $p_n$ is as in \eqref{eq:p ge}. 

Because $p_{n+1}$ is close to $p_n$ for large $n$ and in view of the correspondence \eqref{eq:p,rho} between $\rho$ and $p$, we have now confirmed that   
\begin{itemize}
	\item if $n$ is even then every good value of $\rho$ is symmetric-binary-good; 
	\item if $n$ is odd then, for large $n$, nearly every good value of $\rho$ is symmetric-binary-good.  
\end{itemize}

One may also note here that for large $n$ the lower bound $\rho_{n,\min}$ (defined in \eqref{eq:rho ge}) is close to (but less than) $0$, whereas the lower bound $p_{n,\min}$ is close to (but less than) $\frac12$.

\section*{%\color{blue}
References}
% 
%\begin{minipage}
\bibliography{C:/Users/ipinelis/Documents/pCloudSync/mtu_pCloud_02-02-17/bib_files/citations04-02-21}
%
%\bibliography{P:/pCloudSync/mtu_pCloud_02-02-17/bib_files/citations04-02-21}
%\bibliography{P:/pCloudSync/mtu_pCloud_02-02-17/bib_files/citations01-09-20}
%\bibliography{P:/pCloudSync/mtu_pCloud_02-02-17/bib_files/citations10.13.18a}
%\bibliography{P:/pCloudSync/mtu_pCloud_02-02-17/bib_files/citations12.13.12}
%\bibliography{C:/Users/ipinelis/Dropbox/mtu/bib_files/citations12.13.12}
%
%\bibliographystyle{monthly-mine}
%\bibliographystyle{monthly}
\bibliographystyle{abbrv}

\end{document}